\documentclass{amsart}

\usepackage[margin=1.4in]{geometry}

\usepackage{type1cm}         
\usepackage{graphicx}        
\usepackage{multicol}        
\usepackage[bottom]{footmisc}

\usepackage{amsfonts, amsmath, amssymb}
\usepackage{xspace}          
\usepackage{ytableau}        
\usepackage{pgf}

\usepackage[draft=true]{minted} 
\setminted{breaklines=true}
\definecolor{lbcolor}{rgb}{0.9,0.9,0.9}
\setminted{bgcolor=lbcolor}
\setminted{fontsize=\footnotesize}

\usepackage{mathtools}
\usepackage{booktabs}        
\usepackage{tikz}            
\usetikzlibrary{shapes.geometric}        
\usepackage{csquotes}        

\usepackage[colorinlistoftodos, bordercolor=orange, backgroundcolor=orange!20, linecolor=orange, textsize=scriptsize]{todonotes} 

\usepackage[noend]{algorithmic}

\usepackage{calc,ifthen,alltt,tabularx,capt-of}
\usepackage[linesnumbered,longend,ruled,vlined]{algorithm2e}

\usepackage[bookmarks=false]{hyperref}

\usepackage{longtable}
\usepackage{makecell}
\usepackage{wrapfig}

\usepackage[capitalise, noabbrev]{cleveref} 
\usepackage{enumitem} 

\newcommand\OSCAR{\texttt{OSCAR}\xspace}
\newcommand\Julia{\texttt{Julia}\xspace}

\newcommand\GAP{\texttt{GAP}\xspace}
\newcommand\CHEVIE{\texttt{CHEVIE}\xspace}

\newcommand\Maple{\texttt{Maple}\xspace}

\newcommand\FF{\mathbb{F}}
\newcommand\ZZ{\mathbb{Z}} 
\newcommand\QQ{\mathbb{Q}} 

\newcommand\RR{\mathbb{R}} 

\DeclareMathOperator{\GL}{GL} 
\DeclareMathOperator{\SL}{SL} 
\DeclareMathOperator{\PGL}{PGL} 
\DeclareMathOperator{\PSL}{PSL} 
\DeclareMathOperator{\Trace}{Tr} 

\DeclareUnicodeCharacter{3C0}{$\pi$}
\DeclareUnicodeCharacter{2124}{\ensuremath{\mathbb{Z}}}
\DeclareUnicodeCharacter{2208}{\ensuremath{\in}}
\DeclareUnicodeCharacter{2264}{\ensuremath{\leq}}
\DeclareUnicodeCharacter{1D456}{\ensuremath{i}}

\definecolor{promptColor}{rgb}{0.0,0.0,0.589}
\definecolor{brkpromptColor}{rgb}{0.589,0.0,0.0}
\definecolor{gapinputColor}{rgb}{0.589,0.0,0.0}
\definecolor{gapoutputColor}{rgb}{0.0,0.0,0.0}
\definecolor{darkgreen}{rgb}{0.05,0.6,0.1}
\definecolor{colrem}{rgb}{0,0.7,0}

\newcommand{\tmfloatcontents}{}
\newlength{\tmfloatwidth}
\newcommand{\tmfloat}[5]{
	\renewcommand{\tmfloatcontents}{#4}
	\setlength{\tmfloatwidth}{\widthof{\tmfloatcontents}+1in}
	\ifthenelse{\equal{#2}{small}}
	{\setlength{\tmfloatwidth}{0.45\linewidth}}
	{\setlength{\tmfloatwidth}{\linewidth}}
	\begin{minipage}[#1]{\tmfloatwidth}
		\begin{center}
			\tmfloatcontents
			\captionof{#3}{#5}
		\end{center}
\end{minipage}}

\newcommand{\fd}{.}

\title{Group theory in OSCAR}
\author{Claus Fieker}
\address{RPTU Kaiserslautern-Landau, Kaiserslautern, Germany.}
\email{claus.fieker@rptu.de}

\author{Max Horn}
\address{RPTU Kaiserslautern-Landau, Kaiserslautern, Germany.}
\email{max.horn@rptu.de}

\usepackage{cite}
 
\begin{document}
\begin{abstract}
\OSCAR \cite{OSCAR} is an innovative new computer algebra system which combines and extends the power of its four cornerstone systems - GAP (group theory), Singular (algebra and algebraic geometry), Polymake (polyhedral geometry), and Antic (number theory). Assuming little familiarity with the subject, we give an introduction to computations in group theory using \OSCAR, as a chapter of the upcoming \OSCAR book \cite{OSCAR-book}.
\end{abstract}

\maketitle

\section{Introduction}
\label{ch:cs-groups}

Groups occur everywhere in mathematics, and are thus also central to computer
algebra. In fact, groups were among the first objects to be investigated
mechanically. The famous Todd-Coxeter algorithm for enumerating cosets
in finitely-presented groups \cite{ToddCoxeter:1936} was implemented in
Manchester in the 1950s \cite{Bandler:1956}
on an electronic computer (as opposed to human computers, to
which the term used to refer),
and elsewhere a few years later (see e.g. \cite{Felsch:1961}; for a
more comprehensive overview of this early history, see \cite{Cannon:1969}).

\medskip

For a beginner, or in fact any non-specialist, groups are difficult in
computer algebra as they come in many different representations depending
  both on the intended computations as well as the origin of the group.
While \OSCAR tries to make good choices automatically, it is still important to know
what some of the options are.
Some standard ways to represent groups include:
  \begin{enumerate}
    \item via transformations of some ``space'', such as
      \begin{itemize}
      \item permutations of a finite set, leading to \emph{permutation groups};
      \item linear maps acting on a vector space, leading to \emph{matrix groups};
      \end{itemize}
    \item via generators and relations, leading to \emph{finitely presented groups};
    \item via a multiplication table.
  \end{enumerate}

\begin{figure}[h]
  \centering
\begin{tikzpicture}[scale=0.7, help lines/.style={blue!30,very thin},
 pics/ngon/.style={code={\tikzset{ngon/.cd,#1}
  \node[regular polygon,regular polygon sides=\pgfkeysvalueof{/tikz/ngon/n},
  minimum size=4cm,draw,ngon/border] (\pgfkeysvalueof{/tikz/ngon/name}){};
   \foreach \X in {1,...,\pgfkeysvalueof{/tikz/ngon/n}}
   {\pgfmathsetmacro{\myfillcolor}{{\LstFillCols}[mod(\X-1,5)]}
    \pgfmathsetmacro{\mydrawcolor}{{\LstDrawCols}[mod(\X-1,5)]}
   \path (\pgfkeysvalueof{/tikz/ngon/name}.center) -- (\pgfkeysvalueof{/tikz/ngon/name}.corner \X)
    node[circle,draw=\mydrawcolor,fill=\myfillcolor,ngon/nodes]{}
    node[pos=\pgfkeysvalueof{/tikz/ngon/label pos},transform shape=false](\X){\X};
    }
    }},flip about/.style={/utils/exec=\pgfmathsetmacro{\posangle}{%
    -1*iseven(\pgfkeysvalueof{/tikz/ngon/n})*180/\pgfkeysvalueof{/tikz/ngon/n}+%
    (#1-1)*360/\pgfkeysvalueof{/tikz/ngon/n}},
    rotate=\posangle,xscale=-1,rotate=-1*\posangle},
    ngon/.cd,n/.initial=5,border/.style={very thick},nodes/.style={very thick},
    name/.initial={ngon},label pos/.initial=1.2,
    angle of/.code 2 args=\pgfmathsetmacro{#2}{%
    -1*iseven(\pgfkeysvalueof{/tikz/ngon/n})*180/\pgfkeysvalueof{/tikz/ngon/n}+%
    (#1-1)*360/\pgfkeysvalueof{/tikz/ngon/n}}]
 \edef\LstFillCols{"gray!50","gray!50","gray!30","gray!30","gray!50"}
 \edef\LstDrawCols{"black","black","black","black","black"}

  \draw[help lines]   (-3.9,-2.9) grid (3.9,2.9);
  \begin{scope}[ngon/n=5,flip about=1,transform shape]
   \pgfkeys{tikz/ngon/angle of={1}{\myangle}}
   \draw[dashed] (\myangle+90:2.9) -- (\myangle+270:2.9); 
   \pgfkeys{tikz/ngon/angle of={2}{\myangle}}
   \draw[dashed] (\myangle+90:3.9) -- (\myangle+270:3.9); 
   \pic{ngon};
  \end{scope}
 \end{tikzpicture}

  \caption{A regular pentagon}
  \label{fig:groupsRegularPolygon}
\end{figure}
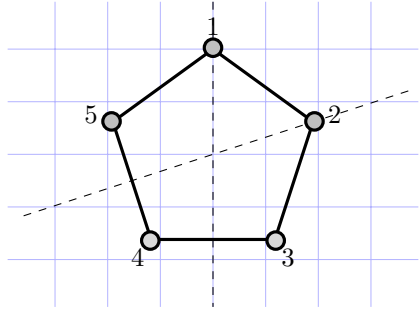

To illustrate these different avenues, consider a regular pentagon in the real plane, centered on the origin and with vertices labeled $1$ to $5$ (see Figure~\ref{fig:groupsRegularPolygon}).
Its symmetries are the reflections $\sigma_i$ at the symmetry axes
through $i$ for $i\in\{1,\dots,5\}$, as well as rotations by multiples of $72$
degrees around the center point. It is well-known that they form a
group under composition,
which has order 10, with 5 reflections and 5 rotations,
and is called the \emph{dihedral group} $D_{10}$.
(We denote the group of symmetries of
a regular $n$-gon as $D_{2n}$, not as $D_n$.)

We now show some different ways of constructing this symmetry group.

\subsection{Permutation groups}

We can interpret the symmetries in purely combinatorial terms by observing
how they permute the vertices of the pentagon. E.g.
$\sigma_2$ leaves vertex $2$ invariant and interchanges $1$ with $3$, and $4$ with $5$.
We can thus describe our group as a \textbf{permutation group}. Here we use that
$\sigma_1$ and $\sigma_2$ already generate it.
\footnote{Note the use of \mintinline{jl}{@permutation_group}
instead of \mintinline{jl}{permutation_group}, which is a \emph{macro} and
allows us to enter the generators conveniently in cycle notation.}

\begin{minted}{jlcon}
julia> G_perm = @permutation_group(5, (2,5)(3,4), (1,3)(4,5))
Permutation group of degree 5
\end{minted}

We can easily verify that it has the correct order and isomorphism type.
\begin{minted}{jlcon}
julia> order(G_perm)
10

julia> describe(G_perm)
"D10"
\end{minted}

The product of two adjacent reflections should give a rotation, and indeed
we can multiply the first and second generator (corresponding to our reflections $\sigma_1$ and $\sigma_2$) to get a $144$-degree clockwise rotation.
Note that in \OSCAR the action is from the right.
\begin{minted}{jlcon}
julia> r = G_perm[1] * G_perm[2]
(1,3,5,2,4)
\end{minted}

Not every permutation of the vertices corresponds to a symmetry of our pentagon.
\begin{minted}{jlcon}
julia> g = perm([1,2,3,5,4])  # 'tabular' notation
(4,5)

julia> g in G_perm
false
\end{minted}
There are several other ways to enter permutations, here are some:
\begin{minted}{jlcon}
julia> G_perm([2,1,5,4,3])  # 'tabular' notation with explicit parent
(1,2)(3,5)

julia> cperm([1,2],[4,5])   # cycle decomposition
(1,2)(4,5)

julia> @perm (1,2)(4,5)     # cycle notation via macro
(1,2)(4,5)
\end{minted}

\subsection{Matrix groups}

A completely different approach is to use the embedding of the pentagon into the real plane $\RR^2$, and then describe
its automorphisms as linear transformations via $2\times2$ matrices.
In Figure~\ref{fig:groupsRegularPolygon} the pentagon can be thought of as centered on the origin $(0,0)$, and we can assign coordinates $(0,1)$ to vertex~$1$.
Then $\sigma_1$ is reflection in the $y$-axis.
To describe a second reflection we would need coordinates for at least one of the other vertices which we don't have. Instead we follow a different approach and exploit that $D_{10}$ is also generated by $\sigma_1$ together with a rotation by 72 degree, i.e., $\frac{2\pi}{5}$ radians.
We can simply write down the corresponding $2\times2$ rotation matrix. To do so we work over an algebraic closure of $\QQ$.
\begin{minted}{jlcon}
julia> K = algebraic_closure(QQ)
Field of algebraic numbers

julia> e = one(K)
Root 1.00000 of x - 1

julia> s, c = sinpi(2*e/5), cospi(2*e/5)
(Root 0.951057 of 16x^4 - 20x^2 + 5, Root 0.309017 of 4x^2 + 2x - 1)

julia> mat_rot = matrix([ c -s ; s c ]);
\end{minted}

\begin{figure}[h]
  \centering
  \includegraphics[height=30mm]{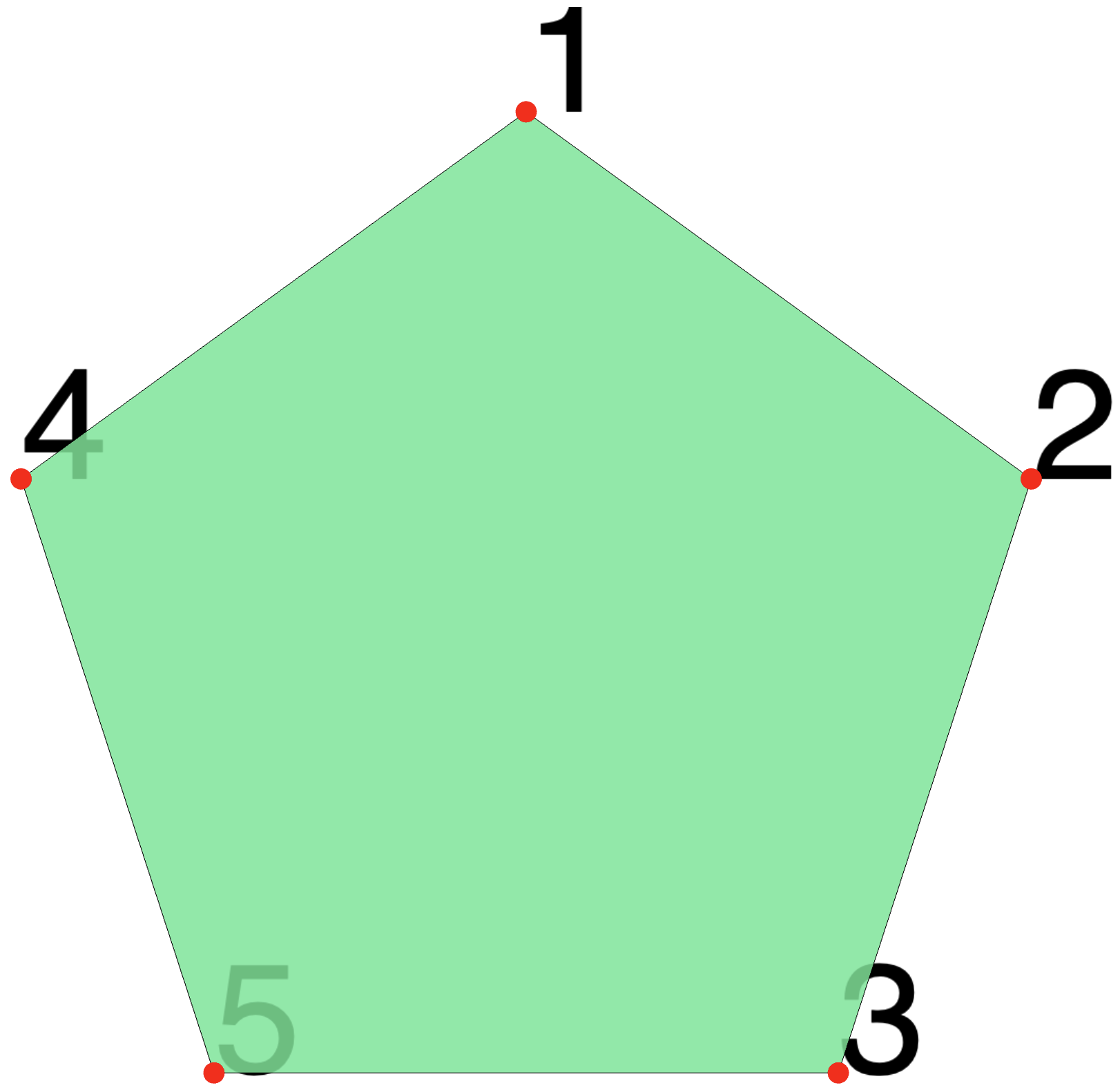}
  \caption{Visualisation of point orbit}
  \label{fig:groupsRegularPolygonScreenshot}
\end{figure}

Together with the matrix for $\sigma_1$ we now have generators for our group:
\begin{minted}{jlcon}
julia> mat_sigma1 = matrix(K, [ -1 0 ; 0 1 ]);

julia> G_mat = matrix_group(mat_rot, mat_sigma1)
Matrix group of degree 2
  over field of algebraic numbers
\end{minted}
This group is isomorphic to the one we constructed
in the previous section:
\begin{minted}{jlcon}
julia> is_isomorphic(G_mat, G_perm)
true
\end{minted}

We can now also compute coordinates for the remaining vertices by computing
the orbit of vertex 1 under the group action.
\begin{minted}{jlcon}
julia> p = K.([0,1]);  # coordinates of vertex 1, expressed over K

julia> orb = orbit(G_mat, *, p)
G-set of
  matrix group of degree 2 over QQBar
  with seeds Vector{QQBarFieldElem}[[Root 0 of x, Root 1.00000 of x - 1]]
\end{minted}
The resulting object is a $G$-set (more on that later),
and at first no actual orbit computation takes place. We can
ask \OSCAR to actually enumerate the orbit by using the standard
\Julia function \mintinline{jl}{collect}.
\begin{minted}{jlcon}
julia> pts = collect(orb)
5-element Vector{Vector{QQBarFieldElem}}:
 [Root 0 of x, Root 1.00000 of x - 1]
 [Root 0.951057 of 16x^4 - 20x^2 + 5, Root 0.309017 of 4x^2 + 2x - 1]
 [Root 0.587785 of 16x^4 - 20x^2 + 5, Root -0.809017 of 4x^2 + 2x - 1]
 [Root -0.951057 of 16x^4 - 20x^2 + 5, Root 0.309017 of 4x^2 + 2x - 1]
 [Root -0.587785 of 16x^4 - 20x^2 + 5, Root -0.809017 of 4x^2 + 2x - 1]
\end{minted}
Often a picture says more than a thousand words (or digits),
so we ask \OSCAR to draw one for us via \mintinline{jl}{visualize(convex_hull(pts))}, with the result in Figure~\ref{fig:groupsRegularPolygonScreenshot}.
(This is actually a 3D visualization, which is why the label $5$ ends up ``behind'' the pentagon.
Note also that the labels in Figure~\ref{fig:groupsRegularPolygonScreenshot}
correspond to the positions in \mintinline{jl}{pts},
which do not coincide with the labels in
Figure~\ref{fig:groupsRegularPolygon}.)
All in all we successfully managed to reconstruct the pentagon we started with.

\subsection{Generic groups}

\OSCAR also supports \textbf{generic groups}, where the user provides
arbitrary objects as generators (or group elements) and specifies a binary function
implementing the multiplication. The result is represented
via a multiplication table -- as such this is only suitable for relatively small
examples.
We again consider our pentagon example, with $K$ the algebraic closure of $\QQ$.

\begin{minted}{jlcon}
julia> R2 = free_module(K, 2) # the "euclidean" plane over K
Vector space of dimension 2 over field of algebraic numbers

julia> A = R2([0,1])
(Root 0 of x, Root 1.00000 of x - 1)

julia> pts = [A*mat_rot^i for i in 0:4];

julia> sigma_1 = hom(R2, R2, [-R2[1], R2[2]])
Module homomorphism
  from vector space of dimension 2 over field of algebraic numbers
  to vector space of dimension 2 over field of algebraic numbers

julia> rot = hom(R2, R2, mat_rot);

julia> G_generic = generic_group(closure([rot, sigma_1], *), *)[1]
Generic group of order 10 with multiplication table
\end{minted}
While generic groups are limited as explained above, in many basic applications
and in teaching, they are perhaps the most natural way of realising a group
on a computer. And in any case it is possible to convert them to other representations,
e.g. into a permutation group.

\begin{minted}{jlcon}
julia> permutation_group(G_generic)
Permutation group of degree 10

julia> is_isomorphic(ans, G_perm)
true
\end{minted}

\subsection{Finitely presented groups}

The symmetry group of the pentagon can also be described as a
\textbf{finitely presented group}, i.e., by  generators and relations.
One way to get such a group is via conversion of an already known group:
\begin{minted}{jlcon}
julia> G_fp = fp_group(G_perm)
Finitely presented group of order 10

julia> gens(G_fp)
2-element Vector{FPGroupElem}:
 F1
 F2

julia> relators(G_fp)
3-element Vector{FPGroupElem}:
 F1^2
 F1^-1*F2*F1*F2^-4
 F2^5
\end{minted}

We may also construct finitely presented groups as quotients of free groups.
Here we use the well-known presentation for $G\cong D_{10}$ as a
reflection group respectively Coxeter group, generated by two
elements of order 2, corresponding to two reflections preserving
our pentagon.
\begin{minted}{jlcon}
julia> F = free_group(2)
Free group of rank 2

julia> G_cox, _ = quo(F, [F[1]^2, F[2]^2, (F[1]*F[2])^5])
(Finitely presented group, Hom: F -> G_cox)

julia> is_isomorphic(G_cox, G_perm)
true
\end{minted}

Small changes in a presentation can have great effects. As an example,
here are two more Coxeter groups whose presentation differs in only one place,
yet they exhibit quite different properties.

\begin{minted}{jlcon}
julia> F = free_group(:a, :b, :c); a,b,c = gens(F);

julia> G, _ = quo(F, [a^2, b^2, c^2, (a*b)^3, (a*c)^2, (b*c)^3])
(Finitely presented group, Hom: F -> G)

julia> H, _ = quo(F, [a^2, b^2, c^2, (a*b)^3, (a*c)^3, (b*c)^3])
(Finitely presented group, Hom: F -> H)

julia> is_finite(G)
true
\end{minted}

The presentations for $G$ and $H$ are quite similar, and \OSCAR easily
determines that $G$ is finite. But trying to do the same for $H$ does
not terminate. In fact $H$ is infinite but the default algorithms are
unable to establish this. So how can we do it?

The problem is that while fp-groups are the ``most natural'' way to present groups in many
applications, it is well known, that in general even to test if a given
fp-group is trivial, let alone if it is finite, is undecidable \cite{zbMATH03115117, zbMATH03141328}.
On the other hand, for many groups arising from applications,
this can be decided quickly.
For mathematicians who are not experts in group theory,
it can often already be helpful to get \emph{some} information about
a given fp-group occurring in their research. The function \mintinline{jl}{describe}
we already saw in use before makes a best effort to produce at least some helpful information about a given group, even though this is mainly impossible in general. Unfortunately
for the second group it still fails to produce anything helpful:
\begin{minted}{jlcon}
julia> describe(G)
"S4"

julia> describe(H)
"a finitely presented group"
\end{minted}

One way to deal with this fundamental problem of fp-groups is to restrict to
working with special kinds of presentation of groups. For example \emph{polycyclic groups}
(pc-groups)
have a specific system of relations that allows for fast computation of
normal forms of words, in particular, to decide if a group element is
trivial or not. On the other hand, polycyclic groups are always
solvable, and hence are not suitable for all applications.

Going back to our group $H$, we can deal with it by finding
a subgroup that is ``obviously'' infinite, namely its second derived subgroup $H''$.
\begin{minted}{jlcon}
julia> H1, _ = derived_subgroup(H)
(Finitely presented group, Hom: H1 -> H)

julia> H2, _ = derived_subgroup(H1)
(Finitely presented group, Hom: H2 -> H1)

julia> describe(H2)
"Z x Z"
\end{minted}

This also establishes that $H$ is actually solvable. The quotient
modulo $H''$ is easily found as
\begin{minted}{jlcon}
julia> describe(quo(H, H2)[1])
"(C3 x C3) : C2"
\end{minted}

\subsection{\ldots and more}

For completeness we mention that of course \OSCAR can also directly produce
a dihedral group if asked for; by default this produces a pc-group but this can
be changed.
\begin{minted}{jlcon}
julia> dihedral_group(10)
Pc group of order 10

julia> dihedral_group(PermGroup, 10)
Permutation group of degree 5
\end{minted}

There are further types of groups available in \OSCAR but for this book
chapter we will consider only the types discussed so far.

While a lot of \OSCAR's group theory capabilities are inherited from \GAP \cite{GAP}
(and indeed anything available in \GAP can be accessed in \OSCAR, see \cref{sect:Outlook}), one
of the compelling advantages of working in \OSCAR is that it is very
convenient to combine them with top-class algorithms for other areas. We will
shortly describe various applications that leverage \OSCAR's number theory
capabilities. Other applications can be seen elsewhere in this book, e.g. in
the chapter on Invariant Theory 
where invariants of permutation
and matrix groups are computed by combining methods from commutative algebra
with group theory tools.

\section{Group actions}

An essential feature of groups is that they \emph{act}
on themselves and on many other kinds of sets and objects. We already
saw this in the form of symmetry groups which act via linear transformations
or vertex permutations on $n$-gons or other geometric objects.
Formally an (right) \emph{action} of a group $G$ on a set $\Omega$ is a function
$f:\Omega\times G\to\Omega,\ (x,g)\mapsto x.g$ such that $x.1_G = x$ and
$(x.g).h = x.(gh)$ for all $x\in\Omega$, $g,h\in G$,
where $1_G$ is the identity element of $G$.

In \OSCAR, actions can be implemented via any \Julia function matching this
signature. However, a pure Julia function does not capture $G$ and
$\Omega$. To deal with this \OSCAR uses \emph{$G$-sets} which are \Julia
objects that can be thought of as \emph{being} the domain $\Omega$,
but stored together with an acting group $G$ and an action function $f$.
Focusing on the domain instead of the action function better matches how one
usually works with and thinks about group actions. On the technical side,
combining everything into a single object allows us to cache properties of the
action and more.

Some examples of $G$-sets in \OSCAR include:
\begin{itemize}
\item the natural $G$-sets of permutation and matrix groups,

\item the set of right cosets of a subgroup $U$ of $G$,

\item conjugacy classes of elements or of subgroups of $G$,

\item vertex or face orbits of the action of a symmetry group on a polyhedron.

\end{itemize}

Here is the $G$-set of the natural action of a permutation group:
\begin{minted}{jlcon}
julia> gset(alternating_group(4))
G-set of
  Alt(4)
  with seeds 1:4
\end{minted}

A slightly more advanced example is the $G$-set of (right) cosets $Ug$, $g\in G$,
of some subgroup $U$ of $G$.
\begin{minted}{jlcon}
julia> G = dihedral_group(6)
Pc group of order 6

julia> U = sub(G, [g for g in gens(G) if order(g) == 2])[1]
Pc group of order 2

julia> r = right_cosets(G, U)
Right cosets of
  pc group of order 2 in
  pc group of order 6
\end{minted}

Note that $r$ is a $G$-set.
We can query it for $G$, its elements (i.e. the underlying set $\Omega$) and the action function $f$.
\begin{minted}{jlcon}
julia> acting_group(r)
Pc group of order 6

julia> collect(r)
3-element Vector{GroupCoset{PcGroup, PcGroupElem}}:
 Right coset of U with representative <identity> of ...
 Right coset of U with representative f2
 Right coset of U with representative f2^2

julia> action_function(r)
* (generic function with 1736 methods)
\end{minted}

As we can see from the last output,
$G$ acts on $r$ via multiplication from the right,
using the normal ($*$)
to denote the operation, hence the large number of methods.
We can now take any element of $G$ and compute how it
permutes the cosets (identified by the indices for which they appear in
the list of right cosets) as follows.
\begin{minted}{jlcon}
julia> permutation(r, G[1])
(2,3)
\end{minted}

This defines an \emph{action homomorphism} $G\to S_{3}$.
If that is all we need then we can also directly ask for it,
without first creating a $G$-set:
\begin{minted}{jlcon}
julia> phi = right_coset_action(G,U)
Group homomorphism
  from pc group of order 6
  to permutation group of degree 3 and order 6

julia> phi(G[1]), phi(G[2])
((2,3), (1,2,3))
\end{minted}

Equipped with this knowledge we can now write a function for
finding a transitive permutation representation of a given group $G$ of minimal degree:
every transitive action of a group on a set $\Omega$ is equivalent to the action of the group on the cosets of a stabilizer $G_a$ for some $a\in\Omega$. We can thus recover all transitive actions by iterating
over all subgroups $U$ of $G$ (up to conjugacy suffices) and checking whether the action of $G$ on the cosets $G/U$ is faithful.
Among these pick one with minimal index
since $[G:U]$ is the degree of the resulting permutation action.
Finally, as an optimization, if $G$ is the symmetric or alternating group of $\Omega$ then we really can't do better and don't bother searching for a permutation
representation of smaller degree.

Combining all that culminates in the following function, which takes a finite group $G$ as input, computes all faithful transitive actions of $G$ and from among those returns one of minimal degree:
\begin{minted}{jlcon}
julia> function optimal_perm_rep(G)
         is_natural_symmetric_group(G) && return hom(G,G,gens(G))
         is_natural_alternating_group(G) && return hom(G,G,gens(G))
         cand = []  # pairs (U,h) with U ≤ G and h a map G -> Sym(G/U)
         for C in subgroup_classes(G)
           U = representative(C)
           h = right_coset_action(G, U)
           is_injective(h) && push!(cand, (U, h))
         end
         return argmax(a -> order(a[1]), cand)[2]
       end;
\end{minted}

We can test it with a small example:
\begin{minted}{jlcon}
julia> U = dihedral_group(8)
Pc group of order 8

julia> optimal_perm_rep(U)
Group homomorphism
  from pc group of order 8
  to permutation group of degree 4 and order 8
\end{minted}

Note that  \mintinline{jl}{permutation_group} and \mintinline{jl}{isomorphism}
also allow conversion of
an arbitrary finite group into an isomorphic permutation group,
but they make no promise to produce something of minimal degree.
\begin{minted}{jlcon}
julia> isomorphism(PermGroup, U)
Group homomorphism
  from pc group of order 8
  to permutation group of degree 8 and order 8

julia> permutation_group(U)
Permutation group of degree 8 and order 8
\end{minted}

\section{Databases}
\label{sect:Databases}

For many research tasks it can be beneficial to have examples at hand.
For this,
libraries of (mainly finite) groups with specific properties are invaluable.
Access to several of these via the cornerstone system \GAP
is directly available in \OSCAR:
\begin{itemize}
\item perfect groups up to order two million,
\cite{Hulpke:2022,PleskenHolt:1989};
\item primitive groups up to degree $8191$,
via the \GAP package \texttt{PrimGrp} \cite{PrimGrp};
\item ``small'' groups of order up to 2000 (except those of order 1024), and many infinite families beyond that, 
via the \GAP packages \texttt{SmallGrp} \cite{SmallGrp}, \texttt{SOTGrps} \cite{SOTGrps},  \texttt{SglPPow} \cite{SglPPow};
\item transitive groups up to degree $48$ (with the caveat that degree 32 and 48
data is huge and has to be downloaded separately; see \cite{trans48} to get an
idea of the tremendous effort it took to compute it),
via the \GAP package \texttt{TransGrp} \cite{TransGrp};
\item
  the Atlas of Group Representations \cite{AGR},
  a collection of permutation and matrix groups
  related to finite simple groups,
  via the \GAP package \texttt{AtlasRep} \cite{AtlasRep};
\item
  the library of character tables from the \texttt{CTblLib} package
  \cite{CTblLib}, which contains among others the character tables from
  the books \cite{CCN85} and \cite{JLPW95}
  (see for example \cref{sect:Characters}).
\end{itemize}

Other {\GAP} data libraries are in principle available in {\OSCAR}
via the {\Julia}-{\GAP} interface (see \cref{sect:Outlook}),
and may get more direct support in the future.
Note that also data libraries in \OSCAR from other areas than group theory can
be interesting in the context of groups, for example by considering the
automorphism groups of objects, see 
the chapter
on Polyhedral Geometry for an example.

All of these databases are the result of hard work by many
researchers over many years. Regretfully, we cannot give a comprehensive
bibliography here. However, the manuals of the referenced packages
include more comprehensive bibliographies on their origins.

Some of these databases are rather large, e.g. the transitive groups of degree
48 take up 30 gigabytes of disk space.
By leveraging the power of \Julia we can
transparently download additional data files on demand from the internet.
\Julia also provides general facilities to manage such ``artifacts'' which we
already use for some of these database. In addition, \Julia provides
convenient interfaces to ``traditional'' database systems (such as PostgreSQL
or MongoDB) which we are actively working on integrating into future \OSCAR
releases to provide more fine grained and flexible access to large
datasets. For example this should eventually make it possible to execute queries such as
``How many transitive groups of degree 48 have point stabilizer isomorphic to
a transitive group'' on a server, then wait a bit and finally get a simple
integer value as result.

As a very small demonstration how one might use these libraries, here
we iterate over all transitive groups of degree $3\leq d\leq 9$ which
are not primitive, and print something for each one that admits a transitive
permutation representation of smaller degree. For this we use the function \mintinline{jl}{optimal_perm_rep} introduced in the previous section.
Here, the notation $(d, n)$ refers to the $n$-transitive group in degree $d$
available via \texttt{transitive\_group(d, n)}.
\begin{minted}{jlcon}
julia> for g in all_transitive_groups(degree => 3:9, !is_primitive)
         h = image(optimal_perm_rep(g))[1]
         if degree(h) < degree(g)
           id = transitive_group_identification(g)
           id_new = transitive_group_identification(h)
           println(id => id_new)
         end
       end
(6, 2) => (3, 2)
(6, 4) => (4, 4)
(6, 7) => (4, 5)
(6, 8) => (4, 5)
(8, 4) => (4, 3)
(8, 13) => (6, 6)
(8, 14) => (4, 5)
(8, 24) => (6, 11)
(9, 4) => (6, 5)
(9, 8) => (6, 9)
\end{minted}

\section{Character theory}
\label{sect:Characters}

For many questions about groups and their representations,
characters are an important tool.
This section shows some applications of complex and Brauer characters.

\subsection{Introduction}

Let $\rho\colon G \to \GL(d, K)$ be a group homomorphism,
where $G$ is finite,
and assume that $K$ is a field of characteristic zero.
Then $\rho$ is called an ordinary (matrix) representation of $G$,
and the \emph{character} of $\rho$ is defined as the map
$$\chi_\rho\colon G \to K, g \mapsto \Trace(\rho(g)).$$

The character $\chi_\rho$ is a much simpler object than $\rho$;
it can be represented by a vector of values in $K$ that is indexed by
the conjugacy classes of elements in $G$.
Still it contains a lot of information about $\rho$:
The degree of $\chi_\rho$ (the value at the identity element of $G$)
equals $d$;
one can compute the characteristic polynomials of matrices $\rho(g)$
from $\chi_\rho$ plus information in which conjugacy classes
the powers of given group elements lie;
the kernel of $\rho$ is the normal subgroup of $G$ that is given by the
union of conjugacy classes of elements $g$ for which $\chi_\rho(g)$
is equal to $d$;
$\rho$ is absolutely irreducible if and only if the norm of $\chi_\rho$
is $1$;
two ordinary representations are equivalent if and only if their
characters are equal.

A group with $n$ conjugacy classes has exactly $n$ absolutely irreducible
representations over a suitably large extension of $K$,
their characters form the \emph{(ordinary) character table} of $G$,
which can be represented by a square matrix of field elements;
by convention,
the rows are the characters, and the columns correspond to the
conjugacy classes of elements in $G$.

To illustrate this we compute the character table of $G = \SL_2(\FF_5)$
(also known as $\SL(2, 5) = \SL_2(5)$ to group theorists).

\begin{minted}{jlcon}
julia> G = SL(2, 5)
SL(2,5)
\end{minted}
\begin{minted}{jlcon}
julia> T = character_table(G)
\end{minted}

We chose to omit the regular output of this command, and instead show in
\cref{fig:character-table} a nicely rendered version by using \LaTeX\ code
produced using
\mintinline{jl}{show(stdout, "text/latex", T)} and some minimal editing.

\begin{figure}[htbp]
\[
\begin{array}{r|rrrrrrrrr}
2 & 3 & 1 & 1 & 3 & 1 & 1 & 1 & 1 & 2 \\
3 & 1 & . & . & 1 & . & . & 1 & 1 & . \\
5 & 1 & 1 & 1 & 1 & 1 & 1 & . & . & . \\
 &  &  &  &  &  &  &  &  &  \\
 & 1a & 10a & 10b & 2a & 5a & 5b & 3a & 6a & 4a \\
2P & 1a & 5b & 5a & 1a & 5b & 5a & 3a & 3a & 2a \\ \hline
\chi_{1} & 1 & 1 & 1 & 1 & 1 & 1 & 1 & 1 & 1 \\
\chi_{2} & 2 & \zeta_{5}^{3} + \zeta_{5}^{2} + 1 & -\zeta_{5}^{3} - \zeta_{5}^{2} & -2 & -\zeta_{5}^{3} - \zeta_{5}^{2} - 1 & \zeta_{5}^{3} + \zeta_{5}^{2} & -1 & 1 & . \\
\chi_{3} & 2 & -\zeta_{5}^{3} - \zeta_{5}^{2} & \zeta_{5}^{3} + \zeta_{5}^{2} + 1 & -2 & \zeta_{5}^{3} + \zeta_{5}^{2} & -\zeta_{5}^{3} - \zeta_{5}^{2} - 1 & -1 & 1 & . \\
\chi_{4} & 3 & -\zeta_{5}^{3} - \zeta_{5}^{2} & \zeta_{5}^{3} + \zeta_{5}^{2} + 1 & 3 & -\zeta_{5}^{3} - \zeta_{5}^{2} & \zeta_{5}^{3} + \zeta_{5}^{2} + 1 & . & . & -1 \\
\chi_{5} & 3 & \zeta_{5}^{3} + \zeta_{5}^{2} + 1 & -\zeta_{5}^{3} - \zeta_{5}^{2} & 3 & \zeta_{5}^{3} + \zeta_{5}^{2} + 1 & -\zeta_{5}^{3} - \zeta_{5}^{2} & . & . & -1 \\
\chi_{6} & 4 & -1 & -1 & 4 & -1 & -1 & 1 & 1 & . \\
\chi_{7} & 4 & 1 & 1 & -4 & -1 & -1 & 1 & -1 & . \\
\chi_{8} & 5 & . & . & 5 & . & . & -1 & -1 & 1 \\
\chi_{9} & 6 & -1 & -1 & -6 & 1 & 1 & . & . & . \\
\end{array}
\]
    \caption{Character table of $\SL_2(\FF_5)$.}
    \label{fig:character-table}
\end{figure}

The rows above the irreducible characters form the header of the table.
It starts with the centralizer orders of the conjugacy classes
in factorized form, each row is labeled by the prime in question.
Then comes a row of conjugacy class names $1a$, $10a$, $10b$, etc.
labeling the columns of the table; each class name consists of the element
order in the class and distinguishing letters.
Finally, the row labeled by $2P$ shows the $2$nd power map, that is,
it contains for each column the name of the class containing the $2$nd
powers of the elements from that class;
more power maps may occur once they get computed.
The character values themselves are shown as sums of roots of unity,
where $\zeta_n$ represents the value $\exp(2 \pi i/n)$;
zeros are abbreviated by dots.
Note that since each character value is a sum of roots of unity,
the entries of the character table lie in a common abelian number field,
in fact $\QQ(\zeta_{|G|})$.

The character table of $G$ contains a lot of information about $G$,
some of which are obvious, others are highly non-trivial:
\begin{itemize}
  \item The number of rows and columns is $n$;
  \item if $G$ is abelian,
nilpotent, supersolvable, solvable, simple, almost simple, perfect;
 \item the lattice of normal subgroups of $G$;
 \item whether the Sylow $p$-subgroups
of $G$ are abelian,
see~\cite{MalleNavarro:2021}.
\end{itemize}

Each character $\chi$ of $G$ is a sum of absolutely irreducible characters,
and once we know the character table of $G$,
we can compute their multiplicities in $\chi$.

The point is that computing the character table of $G$ is usually much easier
than computing the irreducible representations of $G$
(cf. \cref{sect:Brueckner}).
Character-theoretic arguments can be used to avoid dealing with
representations at all, or at least to reduce the work with explicit
representations.

If the characteristic of $K$ is a prime $p$ that divides $|G|$ then the situation is
more involved.
Traces of representing matrices are not sufficient to define characters
in this case.
Instead, one introduces \emph{Brauer characters} that are defined
on the $p$-regular elements of $G$,
and take values in a field of characteristic zero.

\subsection{Character tables related to the Atlas of Finite Groups}

Since a character table can be used to answer many questions about its group,
it makes sense to deal with character tables which do not actually store
a group object -- it is sufficient to know that a given character table
belongs to a particular group, and we can study the group just via its table.

The famous Atlas of Finite Groups \cite{CCN85} contains several hundred
character tables of simple groups and their automorphic and central
extensions,
among them are all such tables for the $26$ sporadic simple groups.
The Atlas defines conventions for the names of the simple groups
and their maximal subgroups,
for the ordering of rows and columns in the character tables, etc.

All character tables from this book are electronically available
in a character table library that can be accessed from {\OSCAR},
with names that fit to the Atlas names,
and with the same ordering of rows and columns as in the book.

This approach has its limitations.
For example, the function \mintinline{jl}{gmodule} for constructing a
representation affording a given character needs a group object associated
to the character table, cf. the example in \cref{sect:SL25}.

In fact, the character table library contains many more character tables
than are shown in the printed Atlas.
For example,
many character tables of maximal subgroups of simple Atlas groups are
available, together with the information how the classes of the subgroup
fuse into the classes of the group;
this allows one to restrict characters to subgroups and to induce characters
from subgroups.

For those library character tables whose groups are contained in the
group libraries described in \cref{sect:Databases},
this information is precomputed and can be used to create a group object
which has the same character table as the given one.

\subsection{Application: Generation questions}

It is known that each finite simple group can be generated by two elements.
A typical question in this context is
from which conjugacy classes one can choose a pair of generators,
and a character theoretic approach is as follows:
The given group $G$ is generated by elements in its conjugacy classes
$C_1$, $C_2$ if and only if not all elements of $C_1 \times C_2$ lie in
maximal subgroups of $G$.
For any subgroup $M$ of $G$ the permutation character $1_M^G$
is nonzero at the class $C$ if and only if $C \cap M$ is empty.
If all maximal subgroups $M$ of $G$ have the property that $1_M^G$ is zero
on $C_1$ or $C_2$ then we see that in fact \emph{each pair} in
$C_1 \times C_2$ generates $G$.

As an example,
let us consider the sporadic simple Janko group $G = J_2$.
The character tables of $G$ and its maximal subgroups are available
in the character table library, and we can compute the permutation
characters $1_M^G$.

\begin{minted}{jlcon}
julia> t = character_table("J2");

julia> mx = character_table.(maxes(t));

julia> pis = [trivial_character(s)^t for s in mx];
\end{minted}

The subgroups in only two classes of maximal subgroups
(of the indices $100$ and $1800$ in $G$) contain elements of order $7$,
and we see that these maximal subgroups do \emph{not} contain elements
of order $5$.

\begin{minted}{jlcon}
julia> ords = orders_class_representatives(t);  println(ords)
[1, 2, 2, 3, 3, 4, 5, 5, 5, 5, 6, 6, 7, 8, 10, 10, 10, 10, 12, 15, 15]

julia> ord7 = findall(is_equal(7), ords);

julia> filt = filter(pi -> pi[ord7[1]] != 0, pis);

julia> degree.(filt)
2-element Vector{QQFieldElem}:
 100
 1800

julia> ord5 = findall(is_equal(5), ords);

julia> s = sum(filt);

julia> all(i -> s[i] == 0, ord5)
true
\end{minted}

Thus we have shown that any pair of elements of the orders $5$ and $7$
in $G$ is a generating set for $G$.

Usually the situation is not so straightforward.
For example, let us consider the question whether $G$ is generated by
a pair of elements of the orders $2$ and $3$ such that their product has
order $7$; groups with this property are called \emph{Hurwitz groups}.

We see that only the maximal subgroups of index $1800$ contain elements
in the conjugacy classes
\texttt{2b} (the second class of elements of order $2$ in $G$)
and \texttt{7a}.
\begin{minted}{jlcon}
julia> (x -> [x["1a"], x["2b"], x["3a"], x["3b"]]).(filt)
2-element Vector{Vector{QQAbElem{AbsSimpleNumFieldElem}}}:
 [100, 0, 10, 4]
 [1800, 20, 0, 6]
\end{minted}

Let $M$ be such a maximal subgroup of index $1800$.
It does not contain elements in the class \texttt{3a},
but this does not help us because no product of a \texttt{2b} and a \texttt{3a}
element can have order $7$.
\begin{minted}{jlcon}
julia> ord3 = findall(is_equal(3), ords);

julia> [class_multiplication_coefficient(t, 3, i, ord7[1]) for i in ord3]
2-element Vector{ZZRingElem}:
 0
 70
\end{minted}

This leaves the class \texttt{3b}.
We see that $G$ contains pairs of \texttt{2b} and \texttt{3b} elements
whose product has order $7$,
and we know that the only maximal subgroups of $G$ that may contain them
are those of index $1800$; these subgroups have the structure $L_3(2).2$,
in the notation defined by the Atlas \cite{CCN85},
as the \mintinline{jl}{identifier} value of their character table shows.%
\footnote{For non-group theorists: $L_3(2).2$ is a group which has a subgroup of index 2 isomorphic to $L_3(2)\cong \GL_3(\FF_2)\cong \SL_3(\FF_2)\cong \PGL_3(\FF_2)\cong \PSL_3(\FF_2)$, in fact $L_3(2).2\cong \PGL_2(\FF_7)$.}
\begin{minted}{jlcon}
julia> s = mx[findfirst(is_equal(filt[2]), pis)];

julia> identifier(s)
"L3(2).2"
\end{minted}

From these subgroups, only class $6$ lies in class \texttt{2b} of $G$,
and only the classes $3$ and $5$ lie in the classes \texttt{3b} and
\texttt{7a}, respectively, of $G$.
It turns out that the subgroups do not contain pairs of \texttt{2b} and
\texttt{3b} elements whose order is $7$.
\begin{minted}{jlcon}
julia> possible_class_fusions(s, t)
1-element Vector{Vector{Int64}}:
 [1, 2, 5, 6, 13, 3, 12, 14, 14]

julia> class_multiplication_coefficient(s, 6, 3, 5)
0
\end{minted}

This proves that $G$ is a Hurwitz group.

It should be noted that character-theoretic computations like the above ones
are not guaranteed to answer the given questions.
Here we were lucky, but otherwise we would have to think about other
strategies.

\section{Applying finitely presented groups towards group extension}

Despite the undecidability of many problems in fp-groups in general, they
are an essential tool in many algorithms.
As an example, given a (permutation) group $G$, via some generators,
a standard problem is to write an arbitrary element of $G$ as a word in the
given generators. This can easily be achieved via fp-groups:
\begin{minted}{jlcon}
julia> G = dihedral_group(PermGroup, 10);

julia> h = epimorphism_from_free_group(G)
Group homomorphism
  from free group of rank 2
  to permutation group of degree 5
\end{minted}
This constructs, as before, the dihedral group with 10 elements, this time
as a permutation group. It is generated by
\begin{minted}{jlcon}
julia> gens(G)
2-element Vector{PermGroupElem}:
 (1,2,3,4,5)
 (2,5)(3,4)
\end{minted}

To decompose an element in $G$ into a word in those generators, we can compute a preimage
under the epimorphism $h$ (this is in general not unique):
\begin{minted}{jlcon}
julia> g = rand(G)
(1,4)(2,3)

julia> w = preimage(h, g)
x1^-1*x2^-1*x1^-3

julia> map_word(w, gens(G))  # evaluate w at generators to get back g
(1,4)(2,3)
\end{minted}

Furthermore, the same map $h$ can also be used to find all relations between
the generators of $G$:
\begin{minted}{jlcon}
julia> k, mk = kernel(h);  # mk is a map from k to domain(h)

julia> ngens(k)
11

julia> mk(k[2])            # pick any kernel generator ...
x1*x2^-2*x1^-1

julia> G[1]*G[-2]^2*G[-1]  # ... and it evaluates to the identity
()
\end{minted}

A second application is building group extensions: given groups $U$ and $V$ any
group extension $X$ satisfying a short exact sequence
$$1 \to U \to X \to V \to 1$$
can naturally be constructed as an fp-group, as a quotient
of a free group of rank \mintinline{jl}{ngens(U) + ngens(V)}.
For split extensions, he relations are then derived from those in $U$, in $V$ and the operation of
$V$ on $U$, in the general case one needs more information.
\begin{minted}{jlcon}
julia> U = cyclic_group(5); V = cyclic_group(2); F = free_group(2);

julia> q = quo(F, [F[1]^5, F[2]^2*F[1]])[1];

julia> describe(q)
"C10"

julia> q, mq = quo(F, [F[1]^5, F[2]^2, F[-1]*F[-2]*F[1]*F[2]*F[1]^2])
(Finitely presented group, Hom: F -> q)

julia> describe(q)
"D10"
\end{minted}

For special classes of fp-groups, typically classes where each group
element has an efficiently computable unique representative (a \emph{normal form}), the
decidability problems vanish. In those classes computations can be
extremely efficient.

Notable examples are groups defined by a confluent rewriting system.
A particularly
nice example for these are \emph{pc-groups}, which are fp-groups given by a special kind of presentation that allows for fast algorithms. On the other hand, pc-groups are
always solvable.\footnote{In fact every pc-group is even polycyclic, a stronger property than solvable. Conversely every polycyclic group is isomorphic to a pc-group.}

To make this concrete, a \emph{finite pc-group} is an fp-group with generators $g_1,\dots,g_n$ and defining relations of the form
\begin{itemize}
  \item for all $1\leq i\leq n$: \emph{power relations} $g_i^{e_i} = $ some word in generators with indices greater than $i$,
  \item for every pair $1\leq i<j\leq n$: \emph{commutator relations}  $[g_i, g_j] = $ some word in generators with indices greater than $i$.
\end{itemize}
The particular shape leads to normal forms for each group element of the
form $g_1^{f_1}\cdots g_n^{f_n}$ with $0\leq f_i < e_i$. In addition to this obvious advantage, there are also
special algorithms for many fundamental problems in group theory available,
see e.g. \cite[Chapter~8]{HandbookComputationalGroup:2005}

In order to use this, the obvious question is to decide if a given group
can be converted into a pc-group (or admits a pc-presentation). This is
the domain of solvable quotient, $p$-quotient and nilpotent
quotient algorithms, and also related to the matrix group recognition problem:
given a group, with limited knowledge, transform it into a more useful
presentation.

\medskip

Next we sketch how to exploit groups with ``nice'' presentations (such as pc-groups) to
compute their cohomology.

\section{Computing cohomology}

Let $G$ be a group. Later we will assume it to be finite. In this section we discuss
how to compute the
second cohomology groups $H^2(G, M)$ for any $G$-module $M$.


Note that most implementations heavily restrict which kinds of modules are supported.
In contrast, for the methods presented here we do not assume
too much about the module: we need homomorphisms, kernels, images and direct sums
and products. As such, our methods work easily for
\begin{itemize}
\item finitely presented abelian groups,
\item (fin. dim.) vector spaces in general,
\item but also e.g. modules over group rings.
\end{itemize}

The basic data type here is \mintinline{jl}{GModule}. Such modules are created by giving
a group, a module, and for each generator an automorphism of the module. For some natural
$G$-modules special constructors are provided, for example for matrix representations.
E.g. for a finite group we can ask for all irreducible representations
in characteristic $0$:
\begin{minted}{jlcon}
julia> irreducible_modules(symmetric_group(3))
3-element Vector{GModule}:
 G-module for Sym(3) acting on vector space of dimension 1 over abelian closure of Q
 G-module for Sym(3) acting on vector space of dimension 1 over abelian closure of Q
 G-module for Sym(3) acting on vector space of dimension 2 over abelian closure of Q
\end{minted}

Or we can construct a module by prescribing an explicit action of the group
on some abelian group. E.g. we can construct the permutation module
of the symmetric group $S_{3}$ over $C_2^3$
\begin{minted}{jlcon}
julia> G = symmetric_group(3);

julia> A = abelian_group([2,2,2])
(Z/2)^3

julia> M = gmodule(G, [hom(A, A, permuted(gens(A), g)) for g in gens(G)])
G-module for G acting on A
\end{minted}

Note that a similar construction would work for any abelian group $A$ (with a suitable operation by $G$). Or, in this specific case, we could just as well have used
\mintinline{jl}{vector_space(GF(2), 3)}. Which we use may have implications for the runtime of subsequent computations with that module, but other than that either
works fine.

\medskip

Now that we know some ways to construct actual modules, we turn to cohomology.

\subsection{$H^2$}
  The computation of $H^2(G,M)$ here is based, similar to \cite{DietrichHulpke:2021}
  and \cite{Holt:2006},
  on the interpretation as group extensions: each element in $H^2(G,M)$
  corresponds to an extension
  $$1\to M \to X \to G\to 1$$
  Elements in $X$ can be represented as pairs $(g, m)$ with the group law
  given by a 2-cochain $\sigma: G^2\to M$ as
  $$(g, m)(h, n) := (gh, m^h + n + \sigma(g, h))$$
  The associativity of the group implies that $\sigma$ has to
  satisfy the cochain condition:
  $$\sigma(g, hk) \sigma(h, k)  = \sigma(g, h)^k \sigma(gh, k), \quad\forall g, h, k\in G$$
  This condition can also be used to find cochains directly. However, naively
  specifying a cochain would require choosing and storing images for $|G|^2$ many elements. How can we overcome this?

  Suppose $G$ is an fp-group, $G = \langle g_1, \ldots, g_n\mid r_1, \ldots, r_s\rangle$, and $M$ is generated by $m_1,\dots,m_\ell$. Then
  \[ X = \langle (g_1, 0), \ldots, (g_n, 0), (1, m_1) , \dots, (1, m_\ell)\rangle \]
  so we have a generating set of $X$, consisting of (lifts of) generators
  for $G$ and (images of) generators of $M$.
  The exact sequence implies that any relation in $M$ also holds in $X$,
  and any relation in $G$, evaluated in $X$
  has to yield an element of $M$. We can thus lift the defining relations of
  $M$ and $G$ to defining relations for $X$. The resulting presentation is
  parametrised by elements $\chi$ of $M^\ell$ while each $r_i$
  yields a linear condition that has to be satisfied by $\chi$.
  Conversely every $\chi\in M^\ell$ satisfying the relations
  induces a presentation for a group $X$ satisfying
  an exact sequence
  $$M\to X\to G \to 1$$
  but in general $M\to X$ will not be injective. In order to achieve this,
  further restrictions on $\chi$ are needed. For a general presentation it is not known
  how to obtain those efficiently, but in some special cases this can be done:
  \begin{itemize}
    \item if the representation  is a confluent rewriting system (RWS), then relations from ``overlaps'' are sufficient;
    \item in case of a pc-presentation, the ``consistency'' relations are well known.
  \end{itemize}
  It should be noted that pc-group relations are very close to a confluent rewriting
  system, in fact translating the pc-relations into a (monoid) RWS shows that
  the consistency relations are exactly those
  from the overlaps in the RWS.

  Let $A$, $B$, $C$ be words in the fixed generators of $G$ such that
  in $G$ we have two relations resp. rewriting rules $AB \to U$ and $BC \to V$.
  Then these translate to relations in $X$ by adding suitable \emph{tails}
  $\chi_1,\chi_2\in M$ to produce
  $$(AB, 0) \to (U, \chi_1), \quad\text{and}\quad (BC, 0) \to (V, \chi_2).$$
  Applying these to $w=(ABC,0)$ we get two possible reductions:
\[
w \to (U, \chi_1)(C, 0) \to (U C, \chi_1^C), \quad\text{ or }\quad
w \to (A, 0)(V,\chi_2) \to (A V, \chi_2).
\]
  Since $AV=UC$ holds in $G$, we deduce that $\chi_1^C = \chi_2$.
  Repeating this process for suitable additional relations gives us the missing
  restrictions on the presentation parameter $\chi\in M^\ell$, and thus
  ultimately allows for an effective computation of $H^2(G,M)$.

\section{Finding irreducible representations of solvable groups}
\label{sect:Brueckner}

A central question in group theory is to find linear representations
of (finite) groups. In theory, this is an easy task: all absolutely irreducible
representations occur as composition factors of the regular
representation. This has been used extremely successfully in positive
characteristic,
via the family of so-called MeatAxe algorithms \cite{Parker:1984,Holt:1998}
and their many refinements. The MeatAxe has been implemented several times
and some version is readily available via the \mintinline{jl}{GModule}
interface:
\begin{minted}{jlcon}
julia> G = dihedral_group(10);

julia> M = regular_gmodule(GF(7), G);

julia> C = composition_factors_with_multiplicity(M)
3-element Vector{Any}:
 (G-module for G acting on vector space of dimension 1 over prime field of characteristic 7, 1)
 (G-module for G acting on vector space of dimension 1 over prime field of characteristic 7, 1)
 (G-module for G acting on vector space of dimension 4 over prime field of characteristic 7, 2)

julia> [is_absolutely_irreducible(x[1]) for x in C]
3-element Vector{Bool}:
 1
 1
 0
\end{minted}

We can easily find the splitting field of the third factor, extend scalars, and then
split it:
\begin{minted}{jlcon}
julia> phi = embed(GF(7), splitting_field(C[3][1]))
Morphism of finite fields
  from prime field of characteristic 7
  to finite field of degree 2 and characteristic 7

julia> M = extension_of_scalars(C[3][1], phi)
G-module for G acting on vector space of dimension 4 over finite field of degree 2 and characteristic 7

julia> composition_factors_with_multiplicity(M)
2-element Vector{Any}:
 (G-module for G acting on vector space of dimension 2 over finite field of degree 2 and characteristic 7, 1)
 (G-module for G acting on vector space of dimension 2 over finite field of degree 2 and characteristic 7, 1)
\end{minted}

And we have all absolutely irreducible representations of $G$ in characteristic $7$.

The situation in characteristic $0$ is different: while deciding
(absolute) irreducibility is possible, splitting is hard,
see \cite{Steel:2012}.
However, in the special case of  pc-groups, this can be done in \OSCAR,
following \cite{Bruckner:1998}:
Let $G = \langle g_1, \ldots, g_r\rangle$ be a pc-group, and set
$U_i = \langle g_i, \ldots, g_r\rangle$. Then the pc-group definition
tells us this is a descending series $G=U_1 \unrhd U_2 \unrhd \dots \unrhd U_{r+1} = 1$. Furthermore, it is straightforward
to arrange for $[U_i : U_{i+1}] = e_i$ to
be prime. The core idea is then that a representation: $\phi:U_{i+1}\to \GL(n, K)$
can easily be extended to a representation of $U_i$:
Suppose there is $X\in \GL(n, K)$ such that setting $\phi(g_i)=X$
yields a valid extension of the representation. Then we must have
for all $y\in U_{i+1}$ that
\[  \phi(y^{g_i}) = \phi(y)^{\phi(g_i)} = \phi(y)^X
\quad\iff\quad \phi(y)X = X\phi(y^{g_i}).
\]
Setting $p:=e_i$, by definition $g_i^p \in U_{i+1}$, thus $\phi(g_i^p)=X^p\in\phi(U_{i+1})$. One can now verify that if $d$ is any zero of the polynomial
$t^p-1$, then $\phi(g_{i+1}) := dX$ also extends $\phi$. For distinct roots
$d$ the resulting extended representations are inequivalent.\footnote{In order to get all roots, we may need to replace $K$ by a larger field. But as we see,
cyclotomic extensions suffice.}

If on the other hand no such $X$ exists,
then inducing $\phi$ to $U_i$ yields an irreducible
representation of degree $np$ over $K$.

Iterating this procedure yields a complete set of absolutely irreducible representations
for $G$.
In those cases where induced representations arise, we can compute which
representations from the previous step are conjugate,
and then by discarding duplicates ensure that
the computed representations are pairwise inequivalent.

Note that in general the representations obtained this way are not
defined over fields of minimal degree:
Let $G = S_{4}$ of order $24$. We can describe it as a pc-group with generators $a$, $b$, $c$ and $d$, power relations
$$a^2, b^3, c^2, d^2$$
and commutator relations
$$[b, a] = b,\; [c,a] = cd,\; [d,a] = cd,\; [c,b] = cd,\; [d,b] = c,\; [d,c]=1.$$

We start with the trivial $1$-dimensional representation of  $\langle c, d\rangle$.
\begin{minted}{jlcon}
julia> G = pc_group(symmetric_group(4));

julia> C = abelian_closure(QQ)[1];

julia> F = free_module(C, 1);

julia> s, ms = sub(G, [G[3], G[4]]);   # subgroup generated by c, d

julia> T = trivial_gmodule(s, F);
\end{minted}

Next we try to extend this to the second generator
$b$ of (relative) order $e_2 = 3$. As $\phi$ at this point is trivial,
$X$ can be chosen as $[1]$.
The distinct roots $d$ of $t^3-1$ are precisely $1$, $\zeta_3$ and $\zeta_3^2$.
For this example we choose to extend $\phi$ via $b\mapsto [\zeta_3]$.
\begin{minted}{jlcon}
julia> z = root_of_unity(C, 3);

julia> ss, mss = sub(G, [G[2], G[3], G[4]]);

julia> M = gmodule(ss, [hom(F, F, [z*F[1]]), action(T, s[1]), action(T, s[2])])
G-module for ss acting on F
\end{minted}

Now trying to extend further to $a$, we see that no $X\in\GL(1,K)$ satisfying the
conditions from above exists:
The commutator relations imply $b^a = b^2$, hence such an $X$ would have to satisfy
\[
\phi(b)X = X\phi(b^a) = X\phi(b^2)
\iff [\zeta_3]X = X[\zeta_3^2]
\]
which is impossible.
So instead we apply the induction procedure to obtain a $2$-dimensional representation:
\begin{minted}{jlcon}
julia> FF = free_module(C, 2);

julia> zm = 0*matrix(action(M, one(ss)));

julia> im = [hom(FF, FF, [matrix(action(M, x)) zm; zm matrix(action(M, preimage(mss, mss(x)^G[1])))]) for x in gens(ss)];
\end{minted}
We obtain the new generator:
\begin{minted}{jlcon}
julia> pushfirst!(im, hom(FF, FF, matrix(C, [0 1; 1 0])));

julia> phi = gmodule(G, im);
\end{minted}
However, we also have
\begin{minted}{jlcon}
julia> character(phi)
class_function(character table of G, QQAbElem{AbsSimpleNumFieldElem}[2, 0, -1, 2, 0])
\end{minted}
so the character is rational, and
\begin{minted}{jlcon}
julia> schur_index(ans)
1
\end{minted}
So, it should be possible to realize this representation over $\QQ$ already.
\noindent Indeed, conjugating by
\begin{minted}{jlcon}
julia> T = matrix(C, [C(1) -z; z+1 z]);
\end{minted}
The result is indeed rational, even integral.
\begin{minted}{jlcon}
julia> [matrix(x)^T for x in action(phi)]
4-element Vector{AbstractAlgebra.Generic.MatSpaceElem{QQAbElem{AbsSimpleNumFieldElem}}}:
 [1 0; -1 -1]
 [0 1; -1 -1]
 [1 0; 0 1]
 [1 0; 0 1]
\end{minted}

It is well known that in general minimal degree realisation fields need not be abelian,
let alone a full cyclotomic field, so finding the matrix $T$ will in general involve
enlarging the field first, see the $\SL_2(\FF_5)$ example in the next section.

From a theoretical point of view, the procedure is straight forward, assuming
the representation is absolutely irreducible. In this case, the representation
defines a central simple algebra, where the center is the character
field. Any such algebra, as an element in the Brauer group, can be
represented by a 2-cocycle with values in a suitable number field. The representation can be realized over a subfield
if and only if the restriction of the 2-cocycle to the fixed group of the subfield
is trivial. In this case, explicitly writing the 2-cocycle as a coboundary as the image of
some 1-cochain allows one to make this explicit. Algorithmically, this
is done using the cohomology described earlier.
To give slightly more details: Let $\psi:G\to \GL_n(K)$ be the representation
and $\sigma:K\to K$ an automorphism fixing the character field. Then
$\psi^\sigma$ is equivalent to $\psi$, hence we can find a matrix
$X_\sigma\in \GL_n(K)$ such that
$$X_\sigma \psi = \psi^\sigma X_\sigma.$$
Let $\tau$ be a second such automorphism, then so is $\sigma\tau$ and there are matrices $X_\tau,X_{\sigma\tau}\in \GL_n(K)$ such that
$$X_\tau \psi = \psi^\tau X_\tau
\quad\text{ and }\quad
X_{\sigma\tau} \psi = \psi^{\sigma\tau} X_{\sigma\tau}
.$$
But we also have
\[ X_\sigma^\tau X_\tau \psi
= X_\sigma^\tau \psi^\tau X_\tau
= (X_\sigma \psi)^\tau X_\tau
= (\psi^\sigma X_\sigma)^\tau X_\tau
= \psi^{\sigma\tau} X_\sigma^\tau X_\tau.
\]
Hence $X_\sigma^\tau X_\tau X_{\sigma\tau}^{-1}$ centralizes $\psi$ and
so $X_\sigma^\tau X_\tau = \mu_{\sigma, \tau} X_{\sigma\tau}$
for some scalar $\mu_{\sigma, \tau}\in K$. As it turns out, computing this
for all (pairs of) automorphisms defines the 2-cycle describing the
algebra in the Brauer group.

The techniques from the previous chapter (plus some number theory, see \cite{Fieker:2009})
can be used to find a 1-cochain $\lambda$ such that $\mu_{\sigma, \tau} = \lambda_\sigma^\tau \lambda_\tau$, thus $\lambda$ is then used to obtain a 1-cochain
of matrices via $X_\sigma \lambda_\sigma$. By Hilbert-90, the 1-cochain of
matrices is actually generated by a single element and this element can be
used to obtain the representation over the smaller field.

\section{Example: Constructing a representation for $\SL_2(\FF_5)$}
\label{sect:SL25}

In this section we return to the perfect group $G = \SL_2(\FF_5)$
whose character table we already saw in \cref{sect:Characters}, \cref{fig:character-table}.
By inspecting this character table we see that there is a unique absolutely irreducible
degree $6$ representation of $G$. We now would like to construct it.

The degree $6$ character is in the final row of the character table.
We compute the corresponding representation over the abelian closure of the rationals.

\begin{minted}{jlcon}
julia> R = gmodule(T[end])
G-module for G acting on vector space of dimension 6 over abelian closure of Q
\end{minted}

To make this more manageable (the abelian closure of $\QQ$ is a complicated object,
much more so than a plain number field), rewrite it over some cyclotomic field (which one is used depends on how the representation was constructed).

\begin{minted}{jlcon}
julia> S = gmodule(CyclotomicField, R)
G-module for G acting on vector space of dimension 6 over cyclotomic field of order 5
\end{minted}

The character in question is rational, it has Schur index $2$,
thus the representation can be written over a quadratic field.
But which one?
First we try to write it over subfields of the given field $\QQ(\zeta_5)$,
in particular over a smallest possible subfield:
\begin{minted}{jlcon}
julia> schur_index(T[end])
2

julia> gmodule_minimal_field(S)
G-module for G acting on vector space of dimension 6 over number field of degree 4 over QQ
\end{minted}

We see that the representation cannot be written over any subfield --
even though there has to be a quadratic field that ``works''.
In order to find out which quadratic field works,
we use more theory: the representation also defines an algebra, namely
the subalgebra of $\QQ(\zeta_5)^{6 \times 6}$ generated by the image
of the representation. Since the character (and hence the representation)
is absolutely irreducible, the algebra is central simple with the
center being the character field.
Classically, those algebras are elements of the Brauer group
of the center ($\QQ$). Algebras defined by matrices in a given
extension $K/\QQ$ are said to be split by $K$, they are represented
by an element of the {\em relative} Brauer group, $\mathop{\mathrm{Br}}(K/\QQ)$.
So let's try that:

\begin{minted}{jlcon}
julia> B, mB = relative_brauer_group(base_ring(S), character_field(S));

julia> B
Relative Brauer group for cyclotomic field of order 5 over number field of degree 1 over QQ
\end{minted}

Classically, the (relative) Brauer group has a well known structure:
Let $A$ be central simple over $k$ (split by $K$),
where $k$ and $K$ are number fields, then for each place
$p$ of $K$, we can consider the completion $K_p$ at this place and the
(local) algebra $A \otimes K_p$. This is now a central simple algebra
over $K_p$ and thus an element in the local Brauer group
$\mathop{\mathrm{Br}}(K_p) \cong \QQ/\ZZ$. Any such algebra is uniquely defined
by the local invariant. Finally, by Hasse-Brauer-Noether
$$\mathop{\mathrm{Br}(K/k)} \cong \sum_p \QQ/\ZZ$$
where the sum is running over all places $p$ of $K$; for each
place $p$, the local invariant is used.
Since an element can have only finitely many non-zero local invariants,
we can represent any such element by the places with non-trivial
invariants, and the invariants.

Furthermore, this representation also helps us with the actual problem: the
non-trivial invariants describe all possible splitting fields!

Here, \mintinline{jl}{B(S)} performs those steps: it turns the
representation \mintinline{jl}{S} into
an element of $\mathop{\mathrm{Br}}(\QQ(\zeta_5)/\QQ)$.

\begin{minted}{jlcon}
julia> b = B(S)
Element of relative Brauer group of number field of degree 1 over QQ
  <2, 2> -> 1//2 + Z
  <5, 5> -> 0 + Z
  Complex embedding of number field -> 1//2 + Z
\end{minted}

The output indicates that we have non-trivial invariants at $2$ and the
``infinite'' place, with both values $1/2 + \ZZ$ of order $2$. So any splitting
field must be totally complex and all completions over primes above $2$ must
also have even degree (the local Schur-indices at $2$ and $\infty$ are $2$).

By the theorem of Grunwald and Wang, there are (many) cyclic extensions
of $\QQ$ matching this, so we feed in the Brauer group element
to obtain a splitting field:

\begin{minted}{jlcon}
julia> K = grunwald_wang(b)
Class field defined mod (<40, 360>, InfPlc{AbsSimpleNumField, AbsSimpleNumFieldEmbedding}[Infinite place corresponding to (Complex embedding corresponding to 1.00 of number field)]) of structure Z/2
\end{minted}

As this field is a ray class field, that is a certain parametrization of
an abelian extension of $\QQ$, we have to convert this to an explicit
extension of $\QQ$:
\begin{minted}{jlcon}
julia> F, _ = absolute_simple_field(number_field(K))
(Number field of degree 2 over QQ, Map: F -> non-simple number field)
\end{minted}

But now we have the next problem, we have a splitting field $F$, but the
representation is over a completely different field! So we first
construct the field $L$ generated by $F$ and
the $5$-th cyclotomic field.

\begin{minted}{jlcon}
julia> L, _, F_to_L = compositum(base_ring(S), F)
(Number field of degree 8 over QQ, Map: cyclotomic field of order 5 -> L, Map: F -> L)
\end{minted}

We write our representation first over $L$,
and then over its quadratic subfield $F$.

\begin{minted}{jlcon}
julia> gmodule_over(F_to_L, gmodule(L, S))
G-module for G acting on vector space of dimension 6 over F
\end{minted}

And thus we have achieved our goal:
an absolutely irreducible degree 6 representation of $\SL_2(\FF_5)$
over a number field of minimal degree 2 over $\QQ$.

\section{Generic character tables}

We have already seen what a character table of a single finite group is.
Some finite groups come in families. For example, the matrix groups $\GL_n(\FF_q)$ or $\SL_n(\FF_q)$, for $n>1$ and $q$ a prime power. These groups have many properties in common.
It turns out that for a fixed rank (say $n=2$)
it is possible to parametrize
the conjugacy classes and irreducible characters of these group in terms of $q$,
and to write this down into a so-called \emph{generic character table}.
This was first done by Schur for $\SL_2(\FF_q)$.

Generic character tables were so far exclusive to the \CHEVIE system
(see e.g. \cite{MR1486215} for details),
which involves use of an old version of \Maple
which heavily restricts the access to this data.
Attempts to port this
to e.g. \GAP in the past were never published due to it being considered too slow.

Thanks to \Julia's impressive performance,
we were able to overcome this obstacle with data imported from \CHEVIE.
The code is still in early stages
as we write this, but this situation is rapidly improving -- so be sure to
check \url{https://book.oscar-system.org} for updated examples for this section.

Right now the code is in a separate repository. To use it, enter within a \Julia interactive command-line (REPL):
\begin{minted}{jl}
  julia> import Pkg; Pkg.add("GenericCharacterTables"); using GenericCharacterTables
\end{minted}

In the following we replicate parts of \cite[Section 5.3]{MR1486215},
and strongly recommend to read that concurrently with the present text.

To start we load the generic character table for $\SL_3(q)$ with
$q\not\equiv 1\pmod 3$.
To learn about the origin of the table, we could enter \mintinline{jl}{printinfotab(T)}.
\begin{minted}{jlcon}
julia> T = genchartab("SL3.n1")
Generic character table
  of order q^8 - q^6 - q^5 + q^3
  with 8 irreducible character types
  with 8 class types
  with parameters (a, b, m, n)
\end{minted}

As we can see this table has four parameters in addition to $q$.
The entries of the table are ``generalized cyclotomics'', that is,
linear combinations over $\QQ(q)$ of symbolic ``roots of unity''
depending on the parameters listed above.
\begin{minted}{jlcon}
julia> printval(T,char=4,class=4)
Value of character type 4 on class type
  4: (q + 1) * exp(2π𝑖(1//(q - 1)*a*n)) + (1) * exp(2π𝑖(-2//(q - 1)*a*n))
\end{minted}

Denoting row 4 of the table by $\chi_4$, we note that is not a single character, but a \emph{character type}, describing a whole family of characters.
We can compute their norms, scalar products, and more. For this demonstration
we tensor the second character type $\chi_2$ with itself.
For technical reasons, the result is currently stored
as an additional character type in the table,
and the \mintinline{jl}{tensor!} function returns its index.
Future versions will make this more elegant.
\begin{minted}{jlcon}
julia> h = tensor!(T,2,2)
9
\end{minted}

We may now attempt to decompose this character type $\chi_9$ by computing its scalar product with the irreducible character types.
This returns a ``generic'' scalar product, plus a (possibly empty)
list of parameter exceptions for which the general result may not hold.
For example:
\begin{minted}{jlcon}
julia> scalar(T,4,h)
(0, Set(ParameterException{QQPolyRingElem}[(2*n1)//(q - 1) ∈ ℤ]))
\end{minted}
This scalar product is $0$
except when $q-1$ divides $2n_1$, where $n_1$ indicates the value of the parameter $n$ for the first factor in the scalar product (i.e. the character type $\chi_4$).
For a ``generic decomposition'' we need to compute all the scalar products and
exceptions.
\begin{minted}{jlcon}
julia> print_decomposition(T, h)
Decomposing character 9:
  <1,9> = 1
  <2,9> = 2
  <3,9> = 2
  <4,9> = 0  with possible exceptions:
    (2*n1)//(q - 1) ∈ ℤ
  <5,9> = 0  with possible exceptions:
    (2*n1)//(q - 1) ∈ ℤ
  <6,9> = 0  with possible exceptions:
    (m1 + n1)//(q - 1) ∈ ℤ
    (2*m1 - n1)//(q - 1) ∈ ℤ
    (m1)//(q - 1) ∈ ℤ
    (n1)//(q - 1) ∈ ℤ
    (m1 - n1)//(q - 1) ∈ ℤ
    (m1 - 2*n1)//(q - 1) ∈ ℤ
  <7,9> = 0  with possible exceptions:
    (n1)//(q - 1) ∈ ℤ
  <8,9> = 0  with possible exceptions:
    ((q + 1)*n1)//(q^2 + q + 1) ∈ ℤ
    (q*n1)//(q^2 + q + 1) ∈ ℤ
    (n1)//(q^2 + q + 1) ∈ ℤ
\end{minted}

This suggest a decomposition of $\chi_2\otimes\chi_2$ into $\chi_1+2\chi_2+2\chi_3$ ``in general''.
But comparing the respective degrees, we notice a discrepancy:
\begin{minted}{jlcon}
julia> chardeg(T, lincomb!(T,[1,2,2],[1,2,3]))
2*q^3 + 2*q^2 + 2*q + 1

julia> chardeg(T, h)
q^4 + 2*q^3 + q^2
\end{minted}

To resolve this, we need to work through the exceptions.
Recall that earlier we saw that $\langle\chi_4,\chi_2^2\rangle=0$
except when $q-1$ divides $2n_1$, where $n_1$ was the value of the parameter $n$
used in $\chi_4$. Further restrictions apply to $n$:
\begin{minted}{jlcon}
julia> printcharparam(T,4)
4	n ∈ {1,…, q - 1} except (n)//(q - 1) ∈ ℤ
\end{minted}

So $n$ may take on any value between $1$ and $q-1$ not divisible by $q-1$.
Hence the only possible exception is $n=(q-1)/2$ which can only
occur if $q$ is odd.
It thus makes sense to consider $q$ odd and $q$ even separately.

We demonstrate this for $q$ even. Then $\langle\chi_4,\chi_2^2\rangle=0$
and with a similar argument $\langle\chi_5,\chi_2^2\rangle=0$.
We now construct a copy of the table but with the congruence equation $q\equiv 0\pmod 2$ applied:
\begin{minted}{jlcon}
julia> T2 = setcongruence(T, (0,2));
\end{minted}
Inspecting the list of exceptions for $\langle\chi_6,\chi_2^2\rangle$, the first
occurs when $q-1$ divides $m+n$. To study this case we specialize $m$:
\begin{minted}{jlcon}
julia> (q, (a, b, m, n)) = params(T2);

julia> x = param(T2, "x");  # create an additional "free" variable

julia> speccharparam!(T2, 6, m, -n + (q-1)*x)  # force m = -n (mod q-1)
\end{minted}

Recomputing the scalar product gives a new result
\begin{minted}{jlcon}
julia> s, e = scalar(T2,6,h); s
1

julia> e
Set{ParameterException{QQPolyRingElem}} with 2 elements:
  (2*n1)//(q - 1) ∈ ℤ
  (3*n1)//(q - 1) ∈ ℤ
\end{minted}

The exceptions both cannot occur as $q$ is even and the table we are considering
is only defined for $q\not\equiv 1\pmod 3$.
By working through the other possible exceptions and irreducible character types, and
handling duplicates, one finally obtains
\[
\chi_2^2 = \chi_1+2\chi_2+2\chi_3
    +\frac12\sum_{n=1}^{q-2} \chi_6(n,q-1-n)
    +\frac12\sum_{n=1}^{q} \chi_7(n(q-1)).
\]
Where $\chi_6(n,q-1-n)$ shows that the $6$th character in the table \texttt{T2}
is a family on two parameters: $n$ and $q-1-n$, while $\chi_7$ depends on only
one, namely $n(q-1)$.
A similar result can be obtained for odd $q$ albeit with a few more cases that need
to be dealt with, but all in essentially the same manner.

\section{Outlook}
\label{sect:Outlook}

{\OSCAR} can already now deal with many kinds of groups, and more group
related objects and algorithms are being actively worked on, such as
Coxeter groups, complex reflection groups, Lie algebras, and more.

\medskip

The abilities of \GAP are vast and extended by many packages.
We currently do not wrap all of it nicely in \OSCAR functions.
However more will be made available over time, prioritized by our
needs and the needs of \OSCAR users -- but for that, we need to
know what they are. So please don't hesitate to contact us to let
us know what you are missing, or even contribute your own additions.

In the meantime, you don't need to wait: our
interface between {\Julia} and {\GAP} allows executing arbitrary
{\GAP} code in a {\Julia} session. The interface is bidirectional,
making it very easy to mix and match GAP and \Julia code.
Global {\GAP} variables can be accessed from \Julia via \mintinline{jl}{GAP.Globals},
for example
\begin{minted}{jl}
  julia> GAP.Globals.SmallGroupsInformation(32)
\end{minted}
calls a {\GAP} function that prints an overview of the groups of order $32$,
and entering
\begin{minted}{jl}
  julia> ?GAP.Globals.SmallGroupsInformation
\end{minted}
shows the documentation of this function.

One can also load additional (perhaps not yet distributed)
{\GAP} packages into the {\Julia} session,
via \mintinline{jl}{GAP.Packages.load},
and then use their functionality.
Furthermore files containing {\GAP} code can be loaded via
\mintinline{jl}{GAP.Globals.Read}.

Using \mintinline{jl}{GAP.prompt()}, you can even get a regular \GAP prompt inside the \Julia REPL.
For details about this
and other features, please have a look at the documentation of
the \texttt{GAP.jl} package.

\section*{Acknowledgements}

The authors acknowledge support by the
  German Research Foundation (DFG) -- Project-ID 286237555 --
  within the SFB-TRR 195 ``Symbolic Tools in Mathematics
  and their Applications'', by MaRDI (Mathematical Research
Data Initiative), funded by the Deutsche Forschungsgemeinschaft (DFG), project
number 460135501, NFDI 29/1 ``MaRDI -- Mathematische
Forschungsdateninitiative'' and ``SymbTools'' funded by the state of Rheinland-Pfalz via the Forschungsinitiative.

\bibliographystyle{plain}
\bibliography{main}

\begin{thebibliography}{10}

\bibitem{GAP}
{GAP} {\textendash} {G}roups, {A}lgorithms, and {P}rogramming, {V}ersion
  4.13.0.
\newblock \url{https://www.gap-system.org}, 2024.

\bibitem{Bandler:1956}
P.~A. Bandler.
\newblock {\em A method for enumerating the cosets of an abstract group on a
  digital computer}.
\newblock M.A. thesis, University of Manchester, 1956.

\bibitem{SmallGrp}
Hans~Ulrich Besche, Bettina Eick, and Eamonn O'Brien.
\newblock {SmallGrp}, the gap small groups library, {V}ersion 1.5.3.
\newblock \url{https://gap-packages.github.io/smallgrp/}, 5 2023.
\newblock GAP package.

\bibitem{zbMATH03141328}
William~Werner Boone.
\newblock The word problem.
\newblock {\em Proc. Natl. Acad. Sci. USA}, 44:1061--1065, 1958.

\bibitem{CTblLib}
Thomas Breuer.
\newblock {CTblLib}, the gap character table library, {V}ersion 1.3.7.
\newblock \url{https://www.math.rwth-aachen.de/~Thomas.Breuer/ctbllib}, 1 2024.
\newblock GAP package.

\bibitem{Bruckner:1998}
Herbert Br{\"u}ckner.
\newblock {\em {Algorithmen f{\"u}r endliche aufl{\"o}sbare Gruppen und
  Anwendungen}}.
\newblock PhD thesis, RWTH Aachen, 1998.

\bibitem{Cannon:1969}
John~J. Cannon.
\newblock Computers in group theory: {A} survey.
\newblock {\em Comm. ACM}, 12:3--12, 1969.

\bibitem{CCN85}
J.~H. Conway, R.~T. Curtis, S.~P. Norton, R.~A. Parker, and R.~A. Wilson.
\newblock {\em {$\mathbb{ATLAS}$} of finite groups}.
\newblock Oxford University Press, Eynsham, 1985.
\newblock Maximal subgroups and ordinary characters for simple groups, With
  computational assistance from J. G. Thackray.

\bibitem{OSCAR-book}
Wolfram Decker, Christian Eder, Claus Fieker, Max Horn, and Michael Joswig,
  editors.
\newblock {\em The {C}omputer {A}lgebra {S}ystem {OSCAR}: {A}lgorithms and
  {E}xamples}.
\newblock Algorithms and {C}omputation in {M}athematics. Springer, 2024.

\bibitem{DietrichHulpke:2021}
Heiko Dietrich and Alexander Hulpke.
\newblock Universal covers of finite groups.
\newblock {\em J. Algebra}, 569:681--712, 2021.

\bibitem{Felsch:1961}
H.~Felsch.
\newblock Programmierung der {R}estklassenabz\"{a}hlung einer {G}ruppe nach
  {U}ntergruppen.
\newblock {\em Numer. Math.}, 3:250--256, 1961.

\bibitem{Fieker:2009}
Claus Fieker.
\newblock Minimizing representations over number fields. {II}. {C}omputations
  in the {B}rauer group.
\newblock {\em J. Algebra}, 322(3):752--765, 2009.

\bibitem{MR1486215}
Meinolf Geck, Gerhard Hiss, Frank L\"{u}beck, Gunter Malle, and G\"{o}tz
  Pfeiffer.
\newblock C{HEVIE}---a system for computing and processing generic character
  tables.
\newblock volume~7, pages 175--210. 1996.
\newblock Computational methods in Lie theory (Essen, 1994).

\bibitem{trans48}
Derek Holt, Gordon Royle, and Gareth Tracey.
\newblock The transitive groups of degree 48 and some applications.
\newblock {\em J. Algebra}, 607:372--386, 2022.

\bibitem{Holt:1998}
Derek~F. Holt.
\newblock The {M}eataxe as a tool in computational group theory.
\newblock In {\em The atlas of finite groups: ten years on ({B}irmingham,
  1995)}, volume 249 of {\em London Math. Soc. Lecture Note Ser.}, pages
  74--81. Cambridge Univ. Press, Cambridge, 1998.

\bibitem{Holt:2006}
Derek~F. Holt.
\newblock Cohomology and group extensions in {M}agma.
\newblock In {\em Discovering mathematics with {M}agma}, volume~19 of {\em
  Algorithms Comput. Math.}, pages 221--241. Springer, Berlin, 2006.

\bibitem{HandbookComputationalGroup:2005}
Derek~F Holt, Bettina Eick, and Eamonn~A O'Brien.
\newblock {\em Handbook of Computational Group Theory}.
\newblock Discrete {{Mathematics}} and Its {{Applications}} ({{Boca Raton}}).
  {Chapman and Hall/CRC}, 2005.

\bibitem{PleskenHolt:1989}
Derek~F. Holt and Wilhelm Plesken.
\newblock {\em Perfect groups}.
\newblock Oxford Mathematical Monographs. The Clarendon Press Oxford University
  Press, New York, 1989.

\bibitem{Hulpke:2022}
Alexander Hulpke.
\newblock The perfect groups of order up to two million.
\newblock {\em Math. Comp.}, 91:1007--1017, 2022.

\bibitem{TransGrp}
Alexander Hulpke.
\newblock {TransGrp}, transitive groups library, {V}ersion 3.6.5.
\newblock \url{https://www.math.colostate.edu/~hulpke/transgrp}, 12 2023.
\newblock GAP package.

\bibitem{PrimGrp}
Alexander Hulpke, Colva Roney-Dougal, and Christopher Russell.
\newblock {PrimGrp}, gap primitive permutation groups library, {V}ersion 3.4.4.
\newblock \url{https://gap-packages.github.io/primgrp/}, 2 2023.
\newblock GAP package.

\bibitem{JLPW95}
C.~Jansen, K.~Lux, R.~Parker, and R.~Wilson.
\newblock {\em An atlas of {B}rauer characters}, volume~11 of {\em London
  Mathematical Society Monographs. New Series}.
\newblock The Clarendon Press Oxford University Press, New York, 1995.
\newblock Appendix 2 by T. Breuer and S. Norton, Oxford Science Publications.

\bibitem{MalleNavarro:2021}
Gunter Malle and Gabriel Navarro.
\newblock Brauer's {H}eight {Z}ero {C}onjecture for principal blocks.
\newblock {\em J. Reine Angew. Math.}, 778:119--125, 2021.

\bibitem{zbMATH03115117}
Pyotr~Sergeyevich Novikov.
\newblock {\"U}ber die algorithmische {Unentscheidbarkeit} des {Wortproblems}
  in der {Gruppentheorie}.
\newblock Tr. {Mat}. {Inst}. {Steklova} 44, 140 {S}. (1955)., 1955.

\bibitem{OSCAR}
{OSCAR -- Open Source Computer Algebra Research system}, version 1.0.0, 2024.

\bibitem{SOTGrps}
Eileen Pan.
\newblock {SOTGrps}, constructing and identifying groups of small order type,
  {V}ersion 1.2.
\newblock \url{https://gap-packages.github.io/sotgrps/}, 6 2023.
\newblock GAP package.

\bibitem{Parker:1984}
R.~A. Parker.
\newblock The computer calculation of modular characters (the meat-axe).
\newblock In {\em Computational group theory ({D}urham, 1982)}, pages 267--274.
  Academic Press, London, 1984.

\bibitem{Steel:2012}
Allan Steel.
\newblock {\em Construction of Ordinary Irreducible Representations of Finite
  Groups}.
\newblock PhD thesis, University of Sydney, 1 2012.

\bibitem{ToddCoxeter:1936}
J.~A. Todd and H.~S.~M. Coxeter.
\newblock A practical method for enumerating cosets of a finite abstract group.
\newblock {\em Proceedings of the Edinburgh Mathematical Society},
  5(1):26–34, 1936.

\bibitem{SglPPow}
Michael Vaughan-Lee and Bettina Eick.
\newblock {SglPPow}, database of groups of prime-power order for some
  prime-powers, {V}ersion 2.3.
\newblock \url{https://gap-packages.github.io/sglppow/}, 11 2022.
\newblock GAP package.

\bibitem{AGR}
R.~A. Wilson, P.~Walsh, J.~Tripp, I.~Suleiman, R.~A. Parker, S.~P. Norton,
  S.~Nickerson, S.~Linton, J.~Bray, and R.~Abbott.
\newblock {ATLAS of Finite Group Representations}.
\newblock \url{https://www.atlasrep.org/Atlas/v3}.

\bibitem{AtlasRep}
Robert~A. Wilson, Richard~A. Parker, Simon Nickerson, John~N. Bray, and Thomas
  Breuer.
\newblock {AtlasRep}, a gap interface to the atlas of group representations,
  {V}ersion 2.1.8.
\newblock \url{https://www.math.rwth-aachen.de/~Thomas.Breuer/atlasrep}, 1
  2024.
\newblock GAP package.

\end{thebibliography}

\end{document}